\documentclass[review]{elsarticle}

\usepackage{bm,amssymb,amsmath,mathrsfs}
\usepackage{amsthm}
\usepackage{lineno}
\usepackage[usenames]{color}

\usepackage{dcolumn}

\theoremstyle{plain}

\newtheorem*{theorem*}{Theorem}

\begin{document}

\journal{(internal report CC26-4)}

\begin{frontmatter}

\title{Power series for roots of a trinomial and Kummer-like identities for higher order hypergeometric series}

\author[cc]{S.~R.~Mane}
\ead{srmane001@gmail.com}
\address[cc]{Convergent Computing Inc., P.~O.~Box 561, Shoreham, NY 11786, USA}

\begin{abstract}
We study the trinomial equation $x^n +px +q =0$.
Here $p$ and $q$ are both real and nonzero.
For $n\ge3$, expressions for the roots have been published as hypergeometric series in powers of the parameter $q^{n-1}/p^n$.
For the special case of the cubic ($n=3$), we employ Kummer's identities to derive alternative series solutions in powers of the discriminant $D$,
and also series in powers of $1/D$.
We next derive new series, in powers of $D$ and also in powers of $1/D$, for all $n\ge 3$.
The resulting series suggest identities analogous to Kummer's identities, for higher order hypergeometric series.
\end{abstract}


\end{frontmatter}

\setcounter{equation}{0}
\section{Introduction}\label{sec:intro}
We derive series solutions for the roots of the trinomial equation
\begin{equation}
\label{eq:tri}
x^n +px +q =0 \,.
\end{equation}
Here $p$ and $q$ are both real and nonzero.
The reason for this constraint is simple.
If $p=0$ then $x^n+q=0$ and the roots are the $n^{th}$ roots of $-q$.
If $q=0$ then $x^n+px=0$ and the origin is a root and the remaining roots are the $(n-1)^{th}$ roots of $-p$.
In principle, our series below (as well as those published in \cite{Birkeland_algebraic_1927,Euler_1779,Lambert_1758,Lambert_1770,Mane_16_1})
are also valid for complex nonzero $p$ and $q$, but branch cuts in the complex plane are required to obtained well-defined expressions.
We shall assume $p$ and $q$ are real in this note.
We note the following general symmetries of the roots of eq.~\eqref{eq:tri}.
\begin{enumerate}
\item
Suppose $n$ is even. If we reverse the sign of $p$, eq.~\eqref{eq:tri} is unchanged if we reverse the sign of $x$.
\item
Suppose $n$ is odd. If we reverse the sign of $q$, eq.~\eqref{eq:tri} is unchanged if we reverse the sign of $x$.
\end{enumerate}

Birkeland \cite{Birkeland_algebraic_1927} derived hypergeometric series solutions for the roots of the more general trinomial equation $x^n = gx^s +\beta$.
For the case $n=3$, Birkeland published the explicit series solutions; see \ref{sec:Birk_cubic}.
The more general trinomial equation has also been solved using Fuss-Catalan series, etc.
See for example, early work by Euler \cite{Euler_1779} and Lambert \cite{Lambert_1758,Lambert_1770}.
A summary is given in \cite{Mane_16_1}, which also contains additional references.

Here we derive alternative series solutions of eq.~\eqref{eq:tri} in powers of the discriminant.
The discriminant $D$ of the trinomial in eq.~\eqref{eq:tri} is
\begin{equation}
\label{eq:discr}  
D = (-1)^n(n-1)^{n-1}p^n -n^nq^{n-1} \,.
\end{equation}
We treat the case $n=3$ (the cubic trinomial) in detail in Section \ref{sec:cubic}.
In Sections \ref{sec:dseries1} and \ref{sec:dseries2},
we treat the case $n>3$ and derive series for the roots in powers of $D$ and $1/D$.
Fuss-Catalan and hypergeometric series for the roots are presented in Section \ref{sec:FChyp}.
The Appendices list relevant identities for ${}_2F_1$ hypergeometric functions
and we employ Kummer's identities to transform Birkeland's series solutions for the trinomial cubic \cite{Birkeland_algebraic_1927} to series in powers of $D$ and $1/D$. 
We also derive a different series, in powers of $D$, in \ref{sec:dseries_alt}, for odd $n\ge3$.

We dispense with the quadratic equation $x^2+px+q=0$ here.
The discriminant is $D = p^2-4q$ and the roots $r_\pm$ are given by
\begin{equation}
r_\pm = \frac{-p \pm\sqrt{D}}{2} \,.
\end{equation}
The roots are functions of $\sqrt{D}$, basically a terminating series.
The formalism in this note is unnecessary for a quadratic.
However, note the following. The roots $r_\pm$ can also be written as
\begin{equation}
r_\pm = \frac{2q}{-p \mp\sqrt{D}} \,.
\end{equation}
In this form, the expressions for $r_\pm$ can be expanded as non-terminating series in powers of $\sqrt{D}/p$,
where the series converges absolutely if and only if $|D| < |p|^2$.
However, we do not pursue the matter further.

\setcounter{equation}{0}
\section{Cubic trinomial}\label{sec:cubic}
We process Cardano's solution for the depressed cubic equation $x^3+px+q=0$, which is a trinomial.
The discriminant is $D = -4p^3 -27q^2$.
Define $u$ and $v$ as follows:
\begin{equation}
u = \biggl(-\frac{q}{2} +\frac{\sqrt{-D}}{2\sqrt{27}}\biggr)^{1/3} \,,\qquad
v = \biggl(-\frac{q}{2} -\frac{\sqrt{-D}}{2\sqrt{27}}\biggr)^{1/3} \,.
\end{equation}
Let $\omega = (-1+i\sqrt3)/2$ be a primitive third root of unity. Then the roots are
\begin{equation}
r_0 = u + v \,,\qquad
r_1 = \omega u + \omega^2 v \,,\qquad
r_2 = \omega^2 u + \omega v \,.
\end{equation}
We begin with $r_0$. Use eq.~\eqref{eq:hyp_sum} to obtain
\begin{equation}
\label{eq:r0}  
\begin{split}
  r_0 &= \Bigl(-\frac{q}{2}\Bigr)^{1/3}\biggl[\,\biggl(1 -\frac{\sqrt{-D}}{\sqrt{27}q}\biggr)^{1/3} +\biggl(1 +\frac{\sqrt{-D}}{\sqrt{27}q}\biggr)^{1/3} \,\biggr]
  \\
  &= 2\Bigl(-\frac{q}{2}\Bigr)^{1/3} \,F\Bigl(-\frac16,\frac13;\frac12;-\frac{D}{27q^2}\Bigr) \,.
\end{split}
\end{equation}
We next process the roots $r_1$ and $r_2$.
For definiteness, let us assume $D<0$, so that $\sqrt{-D}$ is real.
(The case $D>0$ follows straightforwardly.)
Then use eqs.~\eqref{eq:hyp_sum} and \eqref{eq:hyp_diff} to obtain
\begin{equation}  
\label{eq:r1}  
\begin{split}
  r_1 &= \Bigl(-\frac{q}{2}\Bigr)^{1/3}\biggl\{ -\frac12 \biggl[ \biggl(1 -\frac{\sqrt{-D}}{\sqrt{27}q}\biggr)^{1/3} +\biggl(1 +\frac{\sqrt{-D}}{\sqrt{27}q}\biggr)^{1/3} \biggr]
  \\
  &\qquad\qquad\qquad
    -i\frac{\sqrt3}{2} \biggl[ -\biggl(1 -\frac{\sqrt{-D}}{\sqrt{27}q}\biggr)^{1/3} +\biggl(1 +\frac{\sqrt{-D}}{\sqrt{27}q}\biggr)^{1/3} \biggr] \biggr\}
  \\
  &= -\Bigl(-\frac{q}{2}\Bigr)^{1/3}\biggl\{ F\Bigl(-\frac16,\frac13;\frac12;-\frac{D}{27q^2}\Bigr)
  -\frac{i}{\sqrt3}\sqrt{-\frac{D}{27q^2}} \,F\Bigl(\frac13,\frac56;\frac32;-\frac{D}{27q^2}\Bigr) \biggr\} \,.
\end{split}
\end{equation}
The expression for $r_2$ is obtained from eq.~\eqref{eq:r1} by reversing the sign of $i$:
\begin{equation}  
\label{eq:r2}  
\begin{split}
  r_2 &= -\Bigl(-\frac{q}{2}\Bigr)^{1/3}\biggl\{ F\Bigl(-\frac16,\frac13;\frac12;-\frac{D}{27q^2}\Bigr)
  +\frac{i}{\sqrt3}\sqrt{-\frac{D}{27q^2}} \,F\Bigl(\frac13,\frac56;\frac32;-\frac{D}{27q^2}\Bigr) \biggr\} \,.
\end{split}
\end{equation}
Next let us process eq.~\eqref{eq:r0} using Kummer's transformation in eq.~\eqref{eq:Kum_b}.
We attach a prime for reasons to be explained later.
\begin{equation}
\label{eq:r0prime}
\begin{split}
  r_0^\prime &= (-4q)^{1/3}\Bigl(1+\frac{D}{27q^2}\Bigr)^{-1/3} \,F\Bigl(\frac23,\frac13;\frac12;-\frac{D}{4p^3}\Bigr)
  \\
  &= (-4q)^{1/3}\Bigl(-\frac{27q^2}{4p^3}\Bigr)^{1/3} \,F\Bigl(\frac23,\frac13;\frac12;-\frac{D}{4p^3}\Bigr)
  \\
  &= \frac{3q}{p} \,F\Bigl(\frac13,\frac23;\frac12;-\frac{D}{4p^3}\Bigr) \,.
\end{split}
\end{equation}
This is a real root and equals equals Zucker's expression (\cite{Zucker}, eq.~(2.7)).
Next let us process the second term in eq.~\eqref{eq:r1}, say $F_s$,
using Kummer's transformation in eq.~\eqref{eq:Kum_a}
\begin{equation}  
\label{eq:Fs}
\begin{split}
  F_s &\equiv \frac{i}{\sqrt3}\Bigl(-\frac{q}{2}\Bigr)^{1/3}\sqrt{-\frac{D}{27q^2}} \,F\Bigl(\frac13,\frac56;\frac32;-\frac{D}{27q^2}\Bigr)
  \\
  &= \frac{i}{\sqrt3}\Bigl(-\frac{q}{2}\Bigr)^{1/3}\sqrt{-\frac{D}{27q^2}}\Bigl(1+\frac{D}{27q^2}\Bigr)^{-1/3} \,F\Bigl(\frac13,\frac23;\frac32;-\frac{D}{4p^3}\Bigr)
  \\
  &= \frac{i}{\sqrt3}\Bigl(-\frac{q}{2}\Bigr)^{1/3}\sqrt{-\frac{D}{27q^2}}\Bigl(-\frac{27q^2}{4p^3}\Bigr)^{1/3} \,F\Bigl(\frac13,\frac23;\frac32;-\frac{D}{4p^3}\Bigr)
  \\
  &= \frac{\sqrt{D}}{6p} \,F\Bigl(\frac13,\frac23;\frac32;-\frac{D}{4p^3}\Bigr) \,.
\end{split}
\end{equation}
The roots $r_1$ and $r_2$ are given by
\begin{equation}
\label{eq:r12prime}
r_1^\prime = -\frac{r_0^\prime}{2} + F_s \,,\qquad
r_2^\prime = -\frac{r_0^\prime}{2} - F_s \,.
\end{equation}
Note the following.
\begin{enumerate}
\item
The series in eqs.~\eqref{eq:r0}, \eqref{eq:r1} and \eqref{eq:r2} converge unconditionally if and only if $|D| < |27q^2|$, whereas
the series in eqs.~\eqref{eq:r0prime} and \eqref{eq:Fs} converge unconditionally if and only if $|D| < |4p^3|$.
\item
Hence we attached primes to the second set of roots: the domains of convergence are different,
even though the two sets of roots are the same (when both sets of series converge).
\item
The constraint $|D| < |27q^2|$ is satisfied only if $p<0$ and $-54q^2 < 4p^3 < 0$.
The constraint $|D| < |4p^3|$ is satisfied for all $p<0$.
Neither series converges if $p>0$.
\item
The root $r_0^\prime$ is always real.
\item
If $D>0$, all three roots $r_0^\prime$, $r_1^\prime$ and $r_2^\prime$ are real.
\item
The presence of the square root $\sqrt{D}$ in the expression for $F_s$ is essential.
If $D>0$, then $r_1^\prime$ and $r_2^\prime$ are real and distinct.
If $D<0$, then $r_1^\prime$ and $r_2^\prime$ form a complex conjugate pair.
If $D=0$, then $r_1^\prime$ and $r_2^\prime$ are real and equal, i.e., a repeated real root.
\end{enumerate}
Next we treat the domain $|D| > 27 |q|^2$.
To do so, let us write (setting $(-1)^{1/3} = -1$ in the expression for $v$)
\begin{equation}
u = \biggl(\frac{\sqrt{-D}}{2\sqrt{27}}\biggr)^{1/3} \biggl(1 -\frac{\sqrt{27}q}{\sqrt{-D}} \biggr)^{1/3} \,,\qquad
v = -\biggl(\frac{\sqrt{-D}}{2\sqrt{27}}\biggr)^{1/3} \biggl(1 +\frac{\sqrt{27}q}{\sqrt{-D}} \biggr)^{1/3} \,.
\end{equation}
To avoid confusion with the series for the roots for $|D|<|27q^2|$,
let us denote the roots by $\rho_0$, $\rho_1$ and $\rho_2$, with an obvious notation.
Then $\rho_0$ is given by
\begin{equation}
\label{eq:rho0}
\begin{split}
  \rho_0 &= u+v
  \\
  &= -\biggl(\frac{\sqrt{-D}}{2\sqrt{27}}\biggr)^{1/3}\biggl[
    \biggl(1 +\frac{\sqrt{27}q}{\sqrt{-D}} \biggr)^{1/3} - \biggl(1 -\frac{\sqrt{27}q}{\sqrt{-D}} \biggr)^{1/3} \biggr]
  \\
  &= \frac{2q}{3} \biggl(\frac{\sqrt{-D}}{2\sqrt{27}}\biggr)^{1/3} \biggl(-\frac{27}{D}\biggr)^{1/2} \,F\Bigl(\frac13,\frac56;\frac32;-\frac{27q^2}{D}\Bigr)
  \\
  &= -q \Bigl(-\frac{4}{D}\Bigr)^{1/3}\,F\Bigl(\frac13,\frac56;\frac32;-\frac{27q^2}{D}\Bigr) \,.
\end{split}
\end{equation}
Next process $\rho_1$ as follows:
\begin{equation}
\label{eq:rho1}
\begin{split}
  \rho_1 &= \omega u +\omega^2 v
  \\
  &= \biggl(\frac{\sqrt{-D}}{2\sqrt{27}}\biggr)^{1/3}\biggl\{
  \frac12 \biggl[\biggl(1 +\frac{\sqrt{27}q}{\sqrt{-D}} \biggr)^{1/3} - \biggl(1 -\frac{\sqrt{27}q}{\sqrt{-D}} \biggr)^{1/3} \biggr]
  \\
  &\qquad\qquad\qquad\qquad
    +i\frac{\sqrt3}{2} \biggl[ \biggl(1 +\frac{\sqrt{27}q}{\sqrt{-D}} \biggr)^{1/3} + \biggl(1 -\frac{\sqrt{27}q}{\sqrt{-D}} \biggr)^{1/3} \biggr]
    \biggr\}
  \\
  &= \biggl(\frac{\sqrt{-D}}{2\sqrt{27}}\biggr)^{1/3}\biggl\{
  \frac{q}{3} \sqrt{-\frac{27}{D}}\, F\Bigl(\frac13,\frac56;\frac32;-\frac{27q^2}{D}\Bigr)
  +i\sqrt3 F\Bigl(-\frac16,\frac13;\frac12;-\frac{27q^2}{D}\Bigr) \biggr\} \,.
\end{split}
\end{equation}
The expression for $\rho_2$ follows similarly, by reversing the sign of $i$:
\begin{equation}
\label{eq:rho2}
\begin{split}
  \rho_2 &= \omega^2 u +\omega v
  \\
  &= \biggl(\frac{\sqrt{-D}}{2\sqrt{27}}\biggr)^{1/3}\biggl\{
  \frac{q}{3} \sqrt{-\frac{27}{D}}\, F\Bigl(\frac13,\frac56;\frac32;-\frac{27q^2}{D}\Bigr)
  -i\sqrt3 F\Bigl(-\frac16,\frac13;\frac12;-\frac{27q^2}{D}\Bigr) \biggr\} \,.
\end{split}
\end{equation}
Let us process the series in eq.~\eqref{eq:rho0} using Kummer's transformation in eq.~\eqref{eq:Kum_d}
\begin{equation}
\begin{split}
  F\Bigl(\frac13,\frac56;\frac32;-\frac{27q^2}{D}\Bigr) 
  &= \frac{\Gamma(\frac32)\Gamma(\frac13)}{\Gamma(\frac76)\Gamma(\frac23)}F\Bigl(\frac13,\frac56;\frac23;1+\frac{27q^2}{D}\Bigr)
\\
&\quad +\Bigl(1+\frac{27q^2}{D}\Bigr)^{1/3}\frac{\Gamma(\frac32)\Gamma(-\frac13)}{\Gamma(\frac13)\Gamma(\frac56)}F\Bigl(\frac76,\frac23;\frac43;1+\frac{27q^2}{D}\Bigr)
\\
&= \frac{3}{2^{2/3}}F\Bigl(\frac13,\frac56;\frac23;-\frac{4p^3}{D}\Bigr)
-\frac{3}{2^{4/3}}\Bigl(-\frac{4p^3}{D}\Bigr)^{1/3}F\Bigl(\frac76,\frac23;\frac43;-\frac{4p^3}{D}\Bigr) \,.
\end{split}
\end{equation}
Denote the processed root by $\rho_0^\prime$, then
\begin{equation}
\label{eq:rhoprime}
\begin{split}
  \rho_0^\prime &= -\frac{3q}{2^{2/3}} \Bigl(-\frac{4}{D}\Bigr)^{1/3}\,F\Bigl(\frac13,\frac56;\frac23;-\frac{4p^3}{D}\Bigr)
  +\frac{3q}{2^{4/3}} \Bigl(-\frac{4}{D}\Bigr)^{1/3}\,\Bigl(-\frac{4p^3}{D}\Bigr)^{1/3}F\Bigl(\frac76,\frac23;\frac43;-\frac{4p^3}{D}\Bigr)
  \\
  &= -3q \Bigl(-\frac{1}{D}\Bigr)^{1/3}\,F\Bigl(\frac13,\frac56;\frac23;-\frac{4p^3}{D}\Bigr)
  +3pq \Bigl(-\frac{1}{D}\Bigr)^{2/3}F\Bigl(\frac76,\frac23;\frac43;-\frac{4p^3}{D}\Bigr) \,.
\end{split}
\end{equation}
Next process the other hypergeometric function in the expression for $\rho_1$:
\begin{equation}
\begin{split}
F\Bigl(-\frac16,\frac13;\frac12;-\frac{27q^2}{D}\Bigr)
&= \frac{\Gamma(\frac12)\Gamma(\frac13)}{\Gamma(\frac23)\Gamma(\frac16)}F\Bigl(-\frac16,\frac13;\frac23;1+\frac{27q^2}{D}\Bigr)
\\
&\quad +\Bigl(1+\frac{27q^2}{D}\Bigr)^{1/3}\frac{\Gamma(\frac12)\Gamma(-\frac13)}{\Gamma(-\frac16)\Gamma(\frac13)}F\Bigl(\frac23,\frac16;\frac43;1+\frac{27q^2}{D}\Bigr)
\\
&= \frac{1}{2^{2/3}}F\Bigl(-\frac16,\frac13;\frac23;-\frac{4p^3}{D}\Bigr)
+\frac{1}{2^{4/3}}\Bigl(-\frac{4p^3}{D}\Bigr)^{1/3}F\Bigl(\frac23,\frac16;\frac43;-\frac{4p^3}{D}\Bigr) \,.
\end{split}
\end{equation}
Let $F_s^\prime$ denote the second term in the expression for $\rho_1$. Then
\begin{equation}
\label{eq:Fsprime}
\begin{split}
F_s^\prime &= \sqrt3 \biggl(\frac{\sqrt{-D}}{2\sqrt{27}}\biggr)^{1/3}F\Bigl(-\frac16,\frac13;\frac12;-\frac{27q^2}{D}\Bigr)
\\
&= \sqrt{3}\Bigl(\frac{\sqrt{-D}}{2\sqrt{27}}\Bigr)^{1/3}\biggl[
\frac{1}{2^{2/3}}F\Bigl(-\frac16,\frac13;\frac23;-\frac{4p^3}{D}\Bigr)
+\frac{1}{2^{4/3}}\Bigl(-\frac{4p^3}{D}\Bigr)^{1/3}F\Bigl(\frac23,\frac16;\frac43;-\frac{4p^3}{D}\Bigr) \biggr]
\\
&= \frac12\biggl[
(-D)^{\frac16}F\Bigl(-\frac16,\frac13;\frac23;-\frac{4p^3}{D}\Bigr)
+\frac{p}{(-D)^{\frac16}}F\Bigl(\frac23,\frac16;\frac43;-\frac{4p^3}{D}\Bigr) \biggr] \,.
\end{split}
\end{equation}
Then
\begin{equation}
\label{eq:rho12prime}
  \rho_1^\prime = -\frac{\rho_0^\prime}{2} +iF_s^\prime \,,\qquad
  \rho_2^\prime = -\frac{\rho_0^\prime}{2} -iF_s^\prime \,.
\end{equation}
Note the following.
\begin{enumerate}
\item
The series in eqs.~\eqref{eq:rho0}, \eqref{eq:rho1} and \eqref{eq:rho2} converge unconditionally if and only if $|D| > |27q^2|$, whereas
the series in eqs.~\eqref{eq:r0prime} and \eqref{eq:Fsprime} converge unconditionally if and only if $|D| > |4p^3|$.
\item
To satisfy the constraint $|D| > |4p^3|$, we must have $D<0$.
Hence the expression $(-D)^{\frac16}$ in eq.~\eqref{eq:Fsprime} is real.
\item
The constraint $|D| > |4p^3|$ also requires $8p^3>-27q^2$, i.e., $p > -(3/2)|q|^{2/3}$.
\item
The root $\rho_0^\prime$ is always real.
\item
The roots $\rho_1^\prime$ and $\rho_2^\prime$ form a complex conjugate pair.
\end{enumerate}

\setcounter{equation}{0}
\section{Discriminant series I}\label{sec:dseries1}
\subsection{Series in powers of the inverse discriminant}\label{sec:ser_inv_D1}
Let us derive a series in powers of $1/D$.
We assume $n$ is odd, so $D = -(n-1)^{n-1}p^n -n^nq^{n-1}$.
Define $\eta = -n^nq^{n-1}/D$ for use below. 
We try a series of the following form, where $G_1$ is a power series 
\begin{equation}
\label{eq:r_eta}  
r = -\frac{q}{p}\,G_1(\eta) \,.
\end{equation}
Let 
\begin{equation}
\label{eq:G1}  
G_1 = 1 + \sum_{k=1}^\infty \alpha_k\eta^k \,.
\end{equation}
Here the $\alpha_k$ are a set of coefficients to be determined.
Substitute into eq.~\eqref{eq:tri} to obtain
\begin{equation}
-\frac{q^n}{p^n}G_1^n -qG_1 + q = 0 \,.
\end{equation}
Observe that
\begin{equation}
\label{eq:qnpn1}
\begin{split}
-\frac{q^n}{p^n} = -\frac{(n-1)^{n-1}q^n}{(n-1)^{n-1}p^n}
&= \frac{(n-1)^{n-1}q^n}{D +n^nq^{n-1}}
\\
&= q\frac{(n-1)^{n-1}q^{n-1}}{D} \,\frac{1}{1 +n^nq^{n-1}/D}
\\
&= -q\frac{(n-1)^{n-1}}{n^n}\frac{\eta}{1-\eta} \,.
\end{split}
\end{equation}
Also define a multinomial sum as follows.
Define $\tilde{M}_0=\tilde{M}_1=0$ and for $k\ge2$:
\begin{equation}
\label{eq:tilde_M_k_alpha}
\tilde{M}_k = \sum_{j_1+2j_2+\dots+(k-1)j_{k-1}=k}\frac{n!}{j_1!j_2!\dots j_{k-1}!(n-j_1-j_2-\dots-j_{k-1})!} \alpha_1^{,j_1}\alpha_2^{j_2}\dots\alpha_{k-1}^{j_{k-1}} \,.
\end{equation}
Then expand in powers of $\eta$ to obtain
\begin{equation}
\begin{split}
0 &= -\frac{q^n}{p^n}\biggl(1 +\sum_{k=1}^\infty \alpha_k\eta^k\biggr)^n -q\biggl(\sum_{k=1}^\infty \alpha_k\eta^k\biggr)
\\
&= -q\frac{(n-1)^{n-1}}{n^n}\eta(1+\eta+\eta^2+\dots)\biggl[1 +n\alpha_1\eta +(n\alpha_2+\tilde{M}_2)\eta^2 +(n\alpha_3+\tilde{M}_3)\eta^3 +\dots\biggr] 
\\
&\qquad
-q\biggl(\sum_{k=1}^\infty \alpha_k\eta^k\biggr) \,.
\end{split}
\end{equation}
We equate terms in powers of $\eta$.
Then
\begin{equation}
\label{eq:alpha_1}
\alpha_1 = -\frac{(n-1)^{n-1}}{n^n} \,.
\end{equation}
For $k>1$ we obtain a recurrence for $\alpha_k$:
\begin{equation}
\begin{split}
\alpha_k &= -\frac{(n-1)^{n-1}}{n^n}\biggl[1 +n\alpha_1 +(n\alpha_2+\tilde{M}_2) +\dots +(n\alpha_{k-1}+\tilde{M}_{k-1})\biggr] 
\\
&= -\frac{(n-1)^{n-1}}{n^n}\Bigl[1 +n(\alpha_1 +\dots +\alpha_{k-1}) +\tilde{M}_2 +\dots +\tilde{M}_{k-1}\Bigr] 
\\
&= \alpha_{k-1} -\frac{(n-1)^{n-1}}{n^n}(n\alpha_{k-1} +\tilde{M}_{k-1}) 
\\
&= \Bigl(1 -\frac{(n-1)^{n-1}}{n^{n-1}}\Bigr)\alpha_{k-1} -\frac{(n-1)^{n-1}}{n^n}\tilde{M}_{k-1} \,. 
\end{split}
\end{equation}

The above idea also works for even $n$.
For even $n$, the discriminant is
$D = (n-1)^{n-1}p^n -n^nq^{n-1}$.
Recall $\eta = -(n^nq^{n-1})/D$.
Now try a series of the form, where $G_2$ is a power series:
\begin{equation}
r = -\frac{q}{p}\,G_2(\chi) \,.
\end{equation}
Substitute in eq.~\eqref{eq:tri} to obtain
\begin{equation}
\frac{q^n}{p^n}G_2^n -pG_2 + q = 0 \,.
\end{equation}
Observe that
\begin{equation}  
\label{eq:qnpn2}
\begin{split}
\frac{q^n}{p^n} = \frac{(n-1)^{n-1}q^n}{(n-1)^{n-1}p^n}
&= \frac{(n-1)^{n-1}q^n}{D +n^nq^{n-1}}
\\
&= q\frac{(n-1)^{n-1}q^{n-1}}{D} \,\frac{1}{1 +n^nq^{n-1}/D}
\\
&= -q\frac{(n-1)^{n-1}}{n^n}\frac{\eta}{1-\eta} \,.
\end{split}
\end{equation}
This is the same as for odd $n$ in eq.~\eqref{eq:qnpn1}.
\emph{All the rest of the formalism above follows. Nothing beyond this point in the previous derivation requires $n$ to be odd.}

Our series for the root $r$ in eq.~\eqref{eq:r_eta} satisfies the symmetries listed in the Introduction.
(i) If $n$ is even, the discriminant is an even function of $p$. Hence $r$ is an odd function of $p$.
(ii) If $n$ is odd, the discriminant is an even function of $q$. Hence $r$ is an odd function of $q$.

\subsection{Series in powers of the discriminant}\label{sec:ser_D1}
Next let us derive a series in powers of $D$.
We begin with odd $n$.
Recall the discriminant has all negative coefficients
$D = -(n-1)^{n-1}p^n -n^nq^{n-1}$.
Here we define $\chi = -D/(n^nq^{n-1})$. Observe that $\chi=1/\eta$.
We try a series of the form, where $G_3$ is a power series and define
\begin{equation}
\label{eq:r_chi}
r = -\frac{q}{p}\,G_3(\chi) \,.
\end{equation}
Let
\begin{equation}
G_3 = 1 + \sum_{k=1}^\infty \beta_k\chi^k \,.
\end{equation}
Substitute in eq.~\eqref{eq:tri} to obtain
\begin{equation}
-\frac{q^n}{p^n}G_3^n -pG_3 + q = 0 \,.
\end{equation}
Leverage eq.~\eqref{eq:qnpn1} to deduce
\begin{equation}
\label{eq:qnpn3}
-\frac{q^n}{p^n} = -q\frac{(n-1)^{n-1}}{n^n}\frac{1/\chi}{1-1/\chi}
= q\frac{(n-1)^{n-1}}{n^n}\frac{1}{1-\chi} \,.
\end{equation}
Also define a multinomial sum as follows.
Define $\tilde{M}_0^\prime=\tilde{M}_1^\prime=0$ and for $k\ge2$:
\begin{equation}
\tilde{M}_k^\prime = \sum_{j_1+2j_2+\dots+(k-1)j_{k-1}=k}\frac{n!}{j_1!j_2!\dots j_{k-1}!(n-j_1-j_2-\dots-j_{k-1})!} \beta_1^{,j_1}\beta_2^{j_2}\dots\beta_{k-1}^{j_{k-1}} \,.
\end{equation}
Then expand in powers of $\chi$ to obtain
\begin{equation}
\label{eq:notcancel}
\begin{split}
0 &= -\frac{q^n}{p^n}\biggl(1 +\sum_{k=1}^\infty \beta_k\chi^k\biggr)^n -q\biggl(\sum_{k=1}^\infty \beta_k\chi^k\biggr)
\\
&= q\frac{(n-1)^{n-1}}{n^n}(1+\chi+\chi^2+\dots)\biggl[1 +n\beta_1\chi +(n\beta_2+\tilde{M}_2^\prime)\chi^2 +(n\beta_3+\tilde{M}_3^\prime)\chi^3 +\dots\biggr] 
\\
&\qquad
-q\biggl(\sum_{k=1}^\infty \beta_k\chi^k\biggr) \,.
\end{split}
\end{equation}
We equate terms in powers of $\chi$.
Then
\begin{equation}
\beta_1 = \Bigl(1-\frac{(n-1)^{n-1}}{n^{n-1}}\Bigr)^{-1} \,.
\end{equation}
For $k>1$ we obtain a recurrence for $\beta_k$:
\begin{equation}
\begin{split}
\beta_k &= \frac{(n-1)^{n-1}}{n^n}\biggl[1 +n\beta_1 +(n\beta_2+\tilde{M}_2^\prime) +\dots +(n\beta_k+\tilde{M}_k^\prime)\biggr] 
\\
&= \frac{(n-1)^{n-1}}{n^n}\Bigl[1 +n(\beta_1 +\dots +\beta_k) +\tilde{M}_2^\prime +\dots +\tilde{M}_k^\prime\Bigr] 
\\
&= \beta_{k-1} +\frac{(n-1)^{n-1}}{n^n}(n\beta_k +\tilde{M}_k^\prime)
\\
\Rightarrow\qquad
\Bigl(1-\frac{(n-1)^{n-1}}{n^{n-1}}\Bigr)\beta_k &= \beta_{k-1} +\frac{(n-1)^{n-1}}{n^n}\tilde{M}_k^\prime
\\
\beta_k &= \Bigl(1-\frac{(n-1)^{n-1}}{n^{n-1}}\Bigr)^{-1}\biggl[\beta_{k-1} +\frac{(n-1)^{n-1}}{n^n}\tilde{M}_k^\prime\biggr] \,.
\end{split}
\end{equation}

\emph{However, the term $q(n-1)^{n-1}/n^n$ in eq.~\eqref{eq:notcancel} was not cancelled in the above procedure.}
The root we derived in eq.~\eqref{eq:r_chi} actually solves the equation
\begin{equation}
r^n + pr +q = \frac{(n-1)^{n-1}}{n^n}q \,.
\end{equation}
\emph{However, all is not lost.}
The equation we have solved is actually
\begin{equation}
r^n + pr +\Bigl(1- \frac{(n-1)^{n-1}}{n^n}\Bigr)q =0 \,.
\end{equation}
Let us define a modified value $q_{\rm eff,1}$ such that
\begin{equation}
\begin{split}
q &= q_{\rm eff,1} \Bigl(1- \frac{(n-1)^{n-1}}{n^n}\Bigr)
\\
\Rightarrow \qquad q_{\rm eff,1} &= q\Bigl(1- \frac{(n-1)^{n-1}}{n^n}\Bigr)^{-1} \,.
\end{split}
\end{equation}
Then employ our series, replacing $q$ by $q_{\rm eff,1}$ in our formalism above.
Call the resulting root $r_{\rm eff,1}$, where
$D_{\rm eff,1} = -(n-1)^{n-1}p^n -n^nq_{\rm eff,1}^{n-1}$ and $\chi_{\rm eff,1} = -D_{\rm eff,1}/(n^nq_{\rm eff,1}^{n-1})$. Then
\begin{equation}
\label{eq:reff1}
r_{\rm eff,1} = -\frac{q_{\rm eff,1}}{p}\,G_3(\chi_{\rm eff,1}) \,.
\end{equation}
Then $r_{\rm eff,1}$ satisfies the original trinomial equation
\begin{equation}
r_{\rm eff,1}^n + pr_{\rm eff,1} +q = 0 \,.
\end{equation}
Admittedly the series for $q_{\rm eff,1}$ is in powers of a modified discriminant, using $q_{\rm eff,1}$ instead of $q$.

The above idea also works for even $n$.
For even $n$, the discriminant is
$D = (n-1)^{n-1}p^n -n^nq^{n-1}$.
Recall $\chi = -D/(n^nq^{n-1})$.
Now try a series of the form, where $G_4$ is a power series:
\begin{equation}
r = -\frac{q}{p}\,G_4(\chi) \,.
\end{equation}
Substitute in eq.~\eqref{eq:tri} to obtain
\begin{equation}
\frac{q^n}{p^n}G_4^n -pG_4 + q = 0 \,.
\end{equation}
Leverage eq.~\eqref{eq:qnpn2} to deduce
\begin{equation}  
\label{eq:qnpn4}
\frac{q^n}{p^n} = -q\frac{(n-1)^{n-1}}{n^n}\frac{1/\chi}{1-1/\chi} = q\frac{(n-1)^{n-1}}{n^n}\frac{1}{1-\chi} \,.
\end{equation}
This is the same as for odd $n$ in eq.~\eqref{eq:qnpn3}.
\emph{All the rest of the formalism above follows. Nothing beyond this point in the previous derivation requires $n$ to be odd.}

Our series for the root $r_{\rm eff,1}$ satisfies the symmetries listed in the Introduction.
(i) If $n$ is even, the discriminant is an even function of $p$. Hence $r_{\rm eff,1}$ is an odd function of $p$.
(ii) If $n$ is odd, the discriminant is an even function of $q$. Hence $r_{\rm eff,1}$ is an odd function of $q$.

\setcounter{equation}{0}
\section{Discriminant series II}\label{sec:dseries2}
\subsection{Series in powers of the discriminant}\label{sec:ser_D2}
The series in Section \ref{sec:dseries1} employed the parameter $\eta = -n^nq^{n-1}/D$ and its inverse $\chi$.
Here we derive series in powers of $\hat{\eta} = (-1)^n(n-1)^{n-1}p^n/D$ and its inverse $\hat{\chi} = D/((-1)^n(n-1)^{n-1}p^n)$.
We begin with the series in $\hat{\chi}$.
We assume $n$ is odd, so $D = -(n-1)^{n-1}p^n -n^nq^{n-1}$ and $\hat{\chi} = -D/((n-1)^{n-1}p^n)$.
We try a series of the following form, where $G_1$ is a power series 
\begin{equation}
\label{eq:r_chi_hat}  
r = -\frac{q}{p}\,G_5(\hat{\chi}) \,.
\end{equation}
Let 
\begin{equation}
G_5 = 1 + \sum_{k=1}^\infty \hat{\beta}_k\hat{\chi}^k \,.
\end{equation}
Here the $\hat{\beta}_k$ are a set of coefficients to be determined.
Substitute into eq.~\eqref{eq:tri} to obtain
\begin{equation}
-\frac{q^n}{p^n}G_5^n -qG_5 + q = 0 \,.
\end{equation}
Observe that
\begin{equation}
\label{eq:qnpn5}
\begin{split}
-\frac{q^n}{p^n} = -q\frac{n^nq^{n-1}}{n^np^n}
  = q\frac{(n-1)^{n-1}p^n+D}{n^np^n}
  = \frac{(n-1)^{n-1}}{n^n} \,q(1-\hat{\chi}) \,.
\end{split}
\end{equation}
Also define a multinomial sum as follows.
Define $\hat{M}_0^\prime=\hat{M}_1^\prime=0$ and for $k\ge2$:
\begin{equation}
\label{eq:Mkhat}  
\hat{M}_k^\prime = \sum_{j_1+2j_2+\dots+(k-1)j_{k-1}=k}\frac{n!}{j_1!j_2!\dots j_{k-1}!(n-j_1-j_2-\dots-j_{k-1})!} \hat{\beta}_1^{,j_1}\hat{\beta}_2^{j_2}\dots\hat{\beta}_{k-1}^{j_{k-1}} \,.
\end{equation}
Then expand in powers of $\hat{\chi}$ to obtain
\begin{equation}
\begin{split}
0 &= -\frac{q^n}{p^n}\biggl(1 +\sum_{k=1}^\infty \hat{\beta}_k\hat{\chi}^k\biggr)^n -q\biggl(\sum_{k=1}^\infty \hat{\beta}_k\hat{\chi}^k\biggr)
\\
&= q\frac{(n-1)^{n-1}}{n^n}(1-\hat{\chi})\biggl[1 +n\hat{\beta}_1\hat{\chi} +(n\hat{\beta}_2+\hat{M}_2^\prime)\hat{\chi}^2 +(n\hat{\beta}_3+\hat{M}_3^\prime)\hat{\chi}^3 +\dots\biggr]
\\
&\qquad
-q\biggl(\sum_{k=1}^\infty \hat{\beta}_k\hat{\chi}^k\biggr) \,.
\end{split}
\end{equation}
We equate terms in powers of $\hat{\chi}$.
Then
\begin{equation}
\hat{\beta}_1 = -\Bigl(1-\frac{(n-1)^{n-1}}{n^{n-1}}\Bigr)^{-1}\frac{(n-1)^{n-1}}{n^{n-1}} \,.
\end{equation}
For $k>1$ we obtain a recurrence for $\hat{\beta}_k$:
\begin{equation}
\begin{split}
\hat{\beta}_k &= \frac{(n-1)^{n-1}}{n^n}\biggl[-(n\hat{\beta}_{k-1}+\hat{M}^\prime_{k-1}) +(n\hat{\beta}_k+\hat{M}^\prime_k)\biggr] 
\\
&= \frac{(n-1)^{n-1}}{n^n}\Bigl[n(\hat{\beta}_k -\hat{\beta}_{k-1}) +\hat{M}^\prime_k-\hat{M}^\prime_{k-1}\Bigr] 
\\
\Rightarrow\qquad
\Bigl(1-\frac{(n-1)^{n-1}}{n^{n-1}}\Bigr)\hat{\beta}_k &= \frac{(n-1)^{n-1}}{n^n} \Bigl[-n\hat{\beta}_{k-1} +\hat{M}^\prime_k-\hat{M}^\prime_{k-1}\Bigr] 
\\
\hat{\beta}_k &= \Bigl(1-\frac{(n-1)^{n-1}}{n^{n-1}}\Bigr)^{-1}\frac{(n-1)^{n-1}}{n^n} \Bigl[-n\hat{\beta}_{k-1} +\hat{M}^\prime_k-\hat{M}^\prime_{k-1}\Bigr] \,.
\end{split}
\end{equation}
As in Section \ref{sec:ser_D1},
there is a constant term $q(n-1)^{n-1}/n^n$ which was not cancelled. 
The root we derived in eq.~\eqref{eq:r_chi_hat} solves the equation
\begin{equation}
r^n + pr +q = \frac{(n-1)^{n-1}}{n^n}q \,.
\end{equation}
The equation we have solved is actually
\begin{equation}
r^n + pr +\Bigl(1- \frac{(n-1)^{n-1}}{n^n}\Bigr)q =0 \,.
\end{equation}
Let us define a modified value $q_{\rm eff,2}$ such that
\begin{equation}
\begin{split}
q &= q_{\rm eff,2} \Bigl(1- \frac{(n-1)^{n-1}}{n^n}\Bigr)
\\
\Rightarrow \qquad q_{\rm eff,2} &= q\Bigl(1- \frac{(n-1)^{n-1}}{n^n}\Bigr)^{-1} \,.
\end{split}
\end{equation}
Then employ our series, replacing $q$ by $q_{\rm eff,2}$ in our formalism above.
Call the resulting root $r_{\rm eff,2}$, where
$D_{\rm eff,2} = -(n-1)^{n-1}p^n -n^nq_{\rm eff,2}^{n-1}$ and $\hat{\chi}_{\rm eff,2} = -D_{\rm eff,2}/((n-1)^{n-1}p^n)$. Then
\begin{equation}
\label{eq:reff2}
r_{\rm eff,2} = -\frac{q_{\rm eff,2}}{p}\,G_6(\chi_{\rm eff,2}) \,.
\end{equation}
Then $r_{\rm eff,2}$ satisfies the original trinomial equation
\begin{equation}
r_{\rm eff,2}^n + pr_{\rm eff,2} +q = 0 \,.
\end{equation}

The above idea also works for even $n$.
For even $n$, the discriminant is
$D = (n-1)^{n-1}p^n -n^nq^{n-1}$ and $\hat{\chi} = D/((n-1)^{n-1}p^n)$.
Try a series of the form, where $G_6$ is a power series:
\begin{equation}
r = -\frac{q}{p}\,G_6(\hat{\chi}) \,.
\end{equation}
Substitute in eq.~\eqref{eq:tri} to obtain
\begin{equation}
\frac{q^n}{p^n}G_6^n -pG_6 + q = 0 \,.
\end{equation}
Observe that
\begin{equation}
\label{eq:qnpn6}
\begin{split}
  \frac{q^n}{p^n} = q\frac{n^nq^{n-1}}{n^np^n}
  = q\frac{(n-1)^{n-1}p^n-D}{n^np^n}
  = \frac{(n-1)^{n-1}}{n^n} \,q(1-\chi) \,.
\end{split}
\end{equation}
The subsequent derivation is then the same as for odd $n$.

Our series for the root $r_{\rm eff,2}$ satisfies the symmetries listed in the Introduction.
(i) If $n$ is even, the discriminant is an even function of $p$. Hence $r_{\rm eff,2}$ is an odd function of $p$.
(ii) If $n$ is odd, the discriminant is an even function of $q$. Hence $r_{\rm eff,2}$ is an odd function of $q$.

\subsection{Series in powers of the inverse discriminant}\label{sec:ser_inv_D2}
Next let us derive a series in powers of $1/D$.
We begin with odd $n$.
Recall the discriminant has all negative coefficients
$D = -(n-1)^{n-1}p^n -n^nq^{n-1}$.
Recall $\hat{\eta} = -(n-1)^{n-1}p^n/D$ for odd $n$ and also that $\hat{\eta} = 1/\hat{\chi}$.
We try a series of the form, where $G_7$ is a power series and define
\begin{equation}
\label{eq:r_hat_eta}
r = -\frac{q}{p}\,G_7(\hat{\eta}) \,.
\end{equation}
Let
\begin{equation}
G_7 = 1 + \sum_{k=1}^\infty \hat{\alpha}_k\hat{\eta}^k \,.
\end{equation}
Here the $\hat{\alpha}_k$ are a set of coefficients to be determined.
Substitute into eq.~\eqref{eq:tri} to obtain
\begin{equation}
-\frac{q^n}{p^n}G_7^n -qG_7 + q = 0 \,.
\end{equation}
Leverage eq.~\eqref{eq:qnpn5} to deduce
\begin{equation}
\label{eq:qnpn7}
-\frac{q^n}{p^n} = \frac{(n-1)^{n-1}}{n^n} \,q\Bigl(1-\frac{1}{\hat{\eta}}\Bigr) \,.
\end{equation}
Also define a multinomial sum as follows.
Define $\hat{M}_0=\hat{M}_1=0$ and for $k\ge2$:
\begin{equation}
\hat{M}_k = \sum_{j_1+2j_2+\dots+(k-1)j_{k-1}=k}\frac{n!}{j_1!j_2!\dots j_{k-1}!(n-j_1-j_2-\dots-j_{k-1})!} \{\hat{\alpha}_1^{,j_1}\hat{\alpha}_2^{j_2}\dots\hat{\alpha}_{k-1}^{j_{k-1}} \,.
\end{equation}
Then expand in powers of $\hat{\eta}$ to obtain
\begin{equation}
\begin{split}
0 &= -\frac{q^n}{p^n}\biggl(1 +\sum_{k=1}^\infty \hat{\alpha}_k\hat{\eta}^k\biggr)^n -q\biggl(\sum_{k=1}^\infty \hat{\alpha}_k\hat{\eta}^k\biggr)
\\
&= q\frac{(n-1)^{n-1}}{n^n}\Bigl(1-\frac{1}{\hat{\eta}}\Bigr)\biggl[1 +n\hat{\alpha}_1\hat{\eta} +(n\hat{\alpha}_2+M_2)\hat{\eta}^2 +(n\hat{\alpha}_3+M_3)\hat{\eta}^3 +\dots\biggr]
\\
&\qquad
-q\biggl(\sum_{k=1}^\infty \hat{\alpha}_k\hat{\eta}^k\biggr) \,.
\end{split}
\end{equation}
There is a leading term in $1/\hat{\eta}$ which does not cancel.
To cancel the constant term (independent of $\hat{\eta}$), we must have
\begin{equation}
\begin{split}
  1 -n\alpha_1 &= 0
  \\
  \Rightarrow\qquad \alpha_1 &= \frac{1}{n} \,.
\end{split}
\end{equation}
For higher powers, we equate terms in powers of $\hat{\eta}$ to obtain a recurrence for $\hat{\alpha}_k$:
\begin{equation}
\begin{split}
\hat{\alpha}_{k-1} &= \frac{(n-1)^{n-1}}{n^n}\biggl[n\hat{\alpha}_{k-1}+\hat{M}_{k-1} -(n\hat{\alpha}_k+\hat{M}_k)\biggr] 
\\
\Rightarrow\qquad
n\hat{\alpha}_k &= \Bigl(1-\frac{n^{n-1}}{(n-1)^{n-1}}\Bigr)n\hat{\alpha}_{k-1} +\hat{M}_{k-1} +\hat{M}_k
\\
\Rightarrow\qquad
\hat{\alpha}_k &= \Bigl(1-\frac{n^{n-1}}{(n-1)^{n-1}}\Bigr)\hat{\alpha}_{k-1} +\frac{\hat{M}_{k-1} +\hat{M}_k}{n} \,.
\end{split}
\end{equation}
Because the leading term in $1/\hat{\eta}$ was not cancelled,
the root in eq.~\eqref{eq:r_hat_eta} solves the equation
\begin{equation}
\begin{split}
r^n +pr +q &= -\frac{q}{\tilde{\eta}}\frac{(n-1)^{n-1}}{n^n} \,.
\\
\Rightarrow \qquad r^n +pr +q\Bigl(1 +\frac{1}{\tilde{\eta}}\frac{(n-1)^{n-1}}{n^n}\Bigr) &= 0 \,.
\end{split}
\end{equation}
Observe that
\begin{equation}
\begin{split}
1 +\frac{1}{\tilde{\eta}}\frac{(n-1)^{n-1}}{n^n} &= 1 -\frac{D}{(n-1)^{n-1}p^n} \frac{(n-1)^{n-1}}{n^n}
\\
&= 1 -\frac{D}{n^np^n}
\\
&= 1 +\frac{(n-1)^{n-1}p^n+n^nq^{n-1}}{n^np^n}
\\
&= 1 +\frac{(n-1)^{n-1}}{n^n}+\frac{q^{n-1}}{p^n} \,.
\end{split}
\end{equation}
Hence we define $q_{\rm eff,3}$ such that
\begin{equation}
q = q_{\rm eff,3} \biggl(1 +\frac{(n-1)^{n-1}}{n^n}+\frac{q_{\rm eff,3}^{n-1}}{p^n}\biggr)
\end{equation}
Then $q_{\rm eff,3}$ must itself satisfy a trinomial equation of degree $n$:
\begin{equation}
q_{\rm eff,3}^n + p^n\Bigl(1 +\frac{(n-1)^{n-1}}{n^n}\Bigr)q_{\rm eff,3} -p^nq = 0 \,.
\end{equation}
Then employ our series, replacing $q$ by $q_{\rm eff,3}$ in our formalism above.
Call the resulting root $r_{\rm eff,3}$, where
$D_{\rm eff,3} = -(n-1)^{n-1}p^n -n^nq_{\rm eff,3}^{n-1}$ and $\hat{\chi}_{\rm eff,3} = -D_{\rm eff,3}/((n-1)^{n-1}p^n)$. Then
\begin{equation}
\label{eq:reff3}
r_{\rm eff,3} = -\frac{q_{\rm eff,3}}{p}\,G_6(\chi_{\rm eff,3}) \,.
\end{equation}
Then $r_{\rm eff,3}$ satisfies the original trinomial equation
\begin{equation}
r_{\rm eff,3}^n + pr_{\rm eff,3} +q = 0 \,.
\end{equation}

The above idea also works for even $n$.
For even $n$, the discriminant is
$D = (n-1)^{n-1}p^n -n^nq^{n-1}$ and $\hat{\eta} = (n-1)^{n-1}p^n/D$.
Try a series of the form, where $G_6$ is a power series:
\begin{equation}
r = -\frac{q}{p}\,G_8(\hat{\eta}) \,.
\end{equation}
Substitute in eq.~\eqref{eq:tri} to obtain
\begin{equation}
\frac{q^n}{p^n}G_8^n -pG_8 + q = 0 \,.
\end{equation}
Observe that
\begin{equation}
\label{eq:qnpn8}
\frac{q^n}{p^n} = q\frac{n^nq^{n-1}}{n^np^n}
= q\frac{(n-1)^{n-1}p^n-D}{n^np^n}
= q\frac{(n-1)^{n-1}}{n^n}\Bigl(1-\frac{1}{\hat{\eta}}\Bigr) \,.
\end{equation}
The subsequent derivation is then the same as for odd $n$.
The equation for $q_{\rm eff,3}$, for all $n$, is
\begin{equation}
q_{\rm eff,3}^n -(-1)^np^n\Bigl(1 +\frac{(n-1)^{n-1}}{n^n}\Bigr)q_{\rm eff,3} +(-1)^np^nq = 0 \,.
\end{equation}

Our series for the root $r_{\rm eff,3}$ satisfies the symmetries listed in the Introduction.
(i) If $n$ is even, the discriminant is an even function of $p$. Hence $r_{\rm eff,3}$ is an odd function of $p$.
(ii) If $n$ is odd, the discriminant is an even function of $q$. Hence $r_{\rm eff,3}$ is an odd function of $q$.

\setcounter{equation}{0}
\section{Fuss-Catalan and hypergeometric series}\label{sec:FChyp}
Birkeland \cite{Birkeland_algebraic_1927} derived hypergeometric series solutions for the roots of the more general trinomial equation $x^n = gx^s +\beta$.
Solutions for the roots of the general trinomial can also be found in \cite{Mane_16_1} (and references therein).
Here we shall confine our attention to the trinomial equation $x^n+px+q=0$ in eq.~\eqref{eq:tri}.
We treat only the simplest case, in which case Fuss-Catalan series for the root is (\cite{Mane_16_1}, eq.~(6.5b))
\begin{equation}
\label{eq:rootFC}
r = -\frac{q}{p}\mathcal{B}\Bigl(n;1;(-1)^n\frac{q^{n-1}}{p^n}\Bigr) \,.
\end{equation}
See \cite{Mane_16_1} for the definition of the Fuss-Catalan series $\mathcal{B}(\dots)$.
The solution for the same root can also be written as a hypergeometric series
\begin{equation}
r = -\frac{q}{p}{}_{n-1}F_{n-2}\Bigl(\frac{1}{n},\frac{2}{n},\dots,\frac{n-1}{n};\frac{2}{n-1},\dots,\frac{n-1}{n},\frac{n}{n-1};(-1)^n\frac{n^n}{(n-1)^{n-1}}\frac{q^{n-1}}{p^n}\Bigr) \,.
\end{equation}
For $n=3$, $4$ and $5$ the series are respectively
\begin{subequations}
\label{eq:roothyp}
\begin{align}
r &= -\frac{q}{p}{}_2F_1\Bigl(\frac13,\frac23;\frac32;-\frac{27q^2}{4p^3}\Bigr) && (n=3)\,,
\\
r &= -\frac{q}{p}{}_3F_2\Bigl(\frac14,\frac24,\frac34;\frac23,\frac43;\frac{256q^3}{27p^4}\Bigr) && (n=4)\,,
\\
r &= -\frac{q}{p}{}_4F_3\Bigl(\frac15,\frac25,\frac35,\frac45;\frac24,\frac34,\frac54;-\frac{3125q^4}{256p^5}\Bigr) && (n=5)\,.
\end{align}
\end{subequations}
By way of example, for $(p,q)=(-5,3.5)$, the series in eqs.~\eqref{eq:reff1}, \eqref{eq:reff2}, \eqref{eq:rootFC}, \eqref{eq:roothyp}
all converge to the same root, which is
\begin{equation}
r \simeq \begin{cases} 0.803908 &\qquad (n=3) \,,\\ 0.770482 &\qquad (n=4) \,,\\ 0.746303 &\qquad (n=5) \,. \end{cases}
\end{equation}
All the series above for the root satisfies the symmetries listed in the Introduction.
(i) If $n$ is even, the root $r$ reverses sign if we reverse the sign $p$.
(ii) If $n$ is odd, the root $r$ reverses sign if we reverse the sign $q$.

As a second example, for $(p,q)=(6,2)$, the series in eqs.~\eqref{eq:r_eta}, \eqref{eq:rootFC}, \eqref{eq:roothyp}
all converge to the same root.
(The root in eq.~\eqref{eq:reff3} was omitted because it requires the solution of an additional trinomial equation to compute the value of $q_{\rm eff,3}$.)
The root is
\begin{equation}
r \simeq \begin{cases} -0.327480 &\qquad (n=3) \,,\\ -0.335444 &\qquad (n=4) \,,\\ -0.332654 &\qquad (n=5) \,. \end{cases}
\end{equation}
All the series above for the root satisfies the symmetries listed in the Introduction.
(i) If $n$ is even, the root $r$ reverses sign if we reverse the sign $p$.
(ii) If $n$ is odd, the root $r$ reverses sign if we reverse the sign $q$.

\newpage
\setcounter{equation}{0}
\section{Comments}\label{sec:comment}
We note the following general observation about series in powers of the discriminant or inverse discriminant.
Consider Kummer's identity in eq.~\eqref{eq:Kum_d}.
It relates $F(a,b;c;z)$ to $F(\dots;1-z)$.
We can interpret $F(\dots;1-z)$ as a Taylor series in $z$ expanded around $z=1$.
Here $z\propto q^{n-1}/p^n$ (for example) and $1-z$ is then proportional to the discriminant $D$.
Consider Kummer's identity in eq.~\eqref{eq:Kum_a} or eq.~\eqref{eq:Kum_b}.
It relates $F(a,b;c;z)$ to $F(\dots;z/(z-1)) = F(\dots;1 + 1/(z-1))$.
Then $1/(z-1) \propto 1/D$.
One can visualize a series in powers of the discriminant $D$ (and also $1/D$)
as a Taylor series in $q^{n-1}/p^n$ expanded around a base point which is not the origin.

Kummer derived identities for ${}_2F_1$ hypergeometric functions, some of which are listed in \ref{sec:hyp}.
The use of those identities in Section \ref{sec:cubic} permitted us to transform a series in $-D/(27q^2)$ to one in $-D/(4p^3)$, etc.
There appear to be no corresponding published identities for higher hypergeometric functions, e.g.~${}_3F_2$ or ${}_4F_3$, etc.
Nevertheless, the series we derived in Sections \ref{sec:dseries1} and \ref{sec:dseries2} do exhibit such relationships, for all $n\ge3$.
Recall $D = (-1)^n(n-1)^{n-1}p^n -n^nq^{n-1}$.
In Section \ref{sec:dseries1}, we derived a series in $\chi = -D/(n^nq^{n-1})$ and in $\eta = -(n^nq^{n-1})/D$.
In Section \ref{sec:dseries2}, we derived a series in $\hat{\chi} = D/((-1)^n(n-1)^{n-1}p^n)$ and in $\hat{\eta} = ((-1)^n(n-1)^{n-1}p^n)/D$.
The series in Section \ref{sec:FChyp} are in powers of $(-1)^nn^nq^{n-1}/((n-1)^{n-1}p^n)$.
Let $z = (-1)^nn^nq^{n-1}/((n-1)^{n-1}p^n)$.
Then
\begin{equation}
\label{eq:z_chi}  
1-z = 1 -\frac{n^nq^{n-1}}{(-1)^n(n-1)^{n-1}p^n} = \frac{(-1)^n(n-1)^{n-1}p^n - n^nq^{n-1}}{(-1)^n(n-1)^{n-1}p^n} = \frac{D}{(-1)^n(n-1)^{n-1}p^n} = \hat{\chi} \,.
\end{equation}
Also
\begin{equation}
\label{eq:z_eta}  
\frac{z}{z-1} = \frac{\displaystyle \frac{n^nq^{n-1}}{(-1)^n(n-1)^{n-1}p^n}}{\displaystyle \frac{n^nq^{n-1}}{(-1)^n(n-1)^{n-1}p^n} -1}
= \frac{n^nq^{n-1}}{n^nq^{n-1} -(-1)^n(n-1)^{n-1}p^n}
= -\frac{n^nq^{n-1}}{D} = \eta \,.
\end{equation}
Next let $z = \chi = -D/(n^nq^{n-1})$. Then
\begin{equation}
  \frac{\chi}{\chi-1} = \frac{\displaystyle -\frac{D}{n^nq^{n-1}}}{\displaystyle -\frac{D}{n^nq^{n-1}}-1}
  = -\frac{D}{-D-n^nq^{n-1}} = \frac{D}{(-1)^n(n-1)^{n-1}p^n} = \hat{\chi} \,.
\end{equation}
Next let $z = \eta = -(n^nq^{n-1})/D$. Then
\begin{equation}
1 -\eta = 1+\frac{n^nq^{n-1}}{D} = \frac{D+n^nq^{n-1}}{D} = \frac{(-1)^n(n-1)^{n-1}p^n}{D} = \hat{\eta} \,.
\end{equation}

We have derived several series, with different arguments, which all satisfy the same trinomial equation $x^n+px+q=0$, for all $n\ge3$.
However, possibly the most significant contribution in this note is not the series themselves but the relations between them.
The evidence suggests that identities analogous to Kummer's do exist, purely in terms of the series themselves, without reference to a specific application.
It is a matter for future research.

\newpage
\setcounter{equation}{0}
\section{Alternative series}\label{sec:alt_series}
Let $r$ denote a root of the original trinomial equation eq.~\eqref{eq:tri}
$x^n+px+q=0$.
We can express this alternatively as
\begin{equation}
\label{eq:r_express}
r = -\frac{1}{p}(q +r^n) \,.
\end{equation}
Recall $\eta = -n^nq^{n-1}/D$.
Consider the root $r$ in eq.~\eqref{eq:r_eta}.
Use also eq.~\eqref{eq:G1} to obtain
\begin{equation}
r = -\frac{q}{p}\,G_1(\eta)
= -\frac{1}{p}\biggl[q + q\alpha_1\eta\Bigl(1 +\frac{\alpha_2}{\alpha_1}\eta +\frac{\alpha_3}{\alpha_1}\eta^2 +\dots\Bigr)\biggr] \,.
\end{equation}
Equating the expressions for $r$ yields
\begin{equation}
r^n = q\alpha_1\,\sum_{k=0}^\infty \frac{\alpha_{k+1}}{\alpha_1}\eta^k \,.
\end{equation}
Hence let us define a new series, say $\tilde{G}_1(\eta)$, such that
\begin{equation}
\tilde{G}_1^n = \sum_{k=0}^\infty \frac{\alpha_{k+1}}{\alpha_1}\eta^k \,.
\end{equation}
By analogy with eq.~\eqref{eq:G1}, let us write
\begin{equation}
\tilde{G}_1 = 1 + \sum_{k=1}^\infty \tilde{\alpha}_k\eta^k \,.
\end{equation}
Then equate terms in powers of $\eta$ to calculate $\tilde{\alpha}_k$ from the $\alpha_j$:
\begin{equation}
G^n = \sum_{k=0}^\infty (n\tilde{\alpha}_k + \tilde{M}_{\rm alt,\,k})\eta^k \equiv \sum_{k=0}^\infty \frac{\alpha_{k+1}}{\alpha_1}\eta^k \,.
\end{equation}
Here $\tilde{M}_{\rm alt,\,k}$ is a multinomial sum in $(\tilde{\alpha}_1,\dots,\tilde{\alpha}_{k-1})$.
It is formally the same sum as $\tilde{M}_k$ in eq.~\eqref{eq:tilde_M_k_alpha},
but with $\tilde{\alpha}_j$ in place of $\alpha_j$.
Hence $\tilde{M}_{\rm alt,\,0}=\tilde{M}_{\rm alt,\,1}=0$ and for $k\ge2$:
\begin{equation}
\tilde{M}_{\rm alt,\,k} = \sum_{j_1+2j_2+\dots+(k-1)j_{k-1}=k}\frac{n!}{j_1!j_2!\dots j_{k-1}!(n-j_1-j_2-\dots-j_{k-1})!} \tilde{\alpha}_1^{,j_1}\tilde{\alpha}_2^{j_2}\dots\tilde{\alpha}_{k-1}^{j_{k-1}} \,.
\end{equation}
Hence $\tilde{\alpha}_0=1$, $\tilde{\alpha}_1=(1/n)(\alpha_2/\alpha_1)$ and for $k\ge2$:
\begin{equation}
\tilde{\alpha}_k = \frac{1}{n}\Bigl(\frac{\alpha_{k+1}}{\alpha_1} -\tilde{M}_k\Bigr) \,.
\end{equation}
Hence the root can be expressed as $r = (q\alpha_1\eta)^{1/n}G_1(\eta)$. This has been numerically confirmed.

Let us relate these results to one of Kummer's identities for hypergeometric series.
Consider eq.~\eqref{eq:Kum_a}, reproduced here for ease of reference:
\begin{equation}
{}_2F_1(a,b;c;z) = (1-z)^{-a}{}_2F_1\Bigl(a,c-b;c;\frac{z}{z-1}\Bigr) \,.
\end{equation}
Consider how this generalizes to higher hypergeometric series.
For our purposes, the identity has the form
\begin{equation}
{}_nF_{n-1}(a,\dots;c,\dots;z) = (1-z)^{-a}{}_nF_{n-1}\Bigl(a,\dots;c,\dots;\frac{z}{z-1}\Bigr) \,.
\end{equation}
Four our specific application to a root of of the trinomial in eq.~\eqref{eq:tri}, we have $a=1/n$.
Also, the argument is $z = (-1)^nn^nq^{n-1}/((n-1)^{n-1}p^n)$.
Then from eq.~\eqref{eq:z_eta}, $z/(z-1) = \eta$, which is the correct argument of the series $G_1$.
Next consider the prefactor $(1-z)^{-a} = (1-z)^{-1/n}$. 
The root of the trinomial is $-(q/p){}_nF_{n-1}(a,\dots;c,\dots;z) = -(q/p)(1-z)^{-a}\dots$.
Using eq.~\eqref{eq:z_chi}, we obtain
\begin{equation}
\label{eq:pref}  
-\frac{q}{p}(1-z)^{-1/n} = -\frac{q}{p}\Bigl(\frac{(-1)^n(n-1)^{n-1}p^n}{D}\Bigr)^{1/n} = q\Bigl(\frac{(n-1)^{n-1}}{D}\Bigr)^{1/n} \,.
\end{equation}
Compare this to our prefactor above (use eq.~\eqref{eq:alpha_1} for $\alpha_1$):
\begin{equation}
(q\alpha_1\eta)^{1/n} = \Bigl(q\,\frac{(n-1)^{n-1}}{n^n}\,\frac{n^nq^{n-1}}{D}\Bigr)^{1/n} = q\Bigl(\frac{(n-1)^{n-1}}{D}\Bigr)^{1/n} \,.
\end{equation}
This exactly matches eq.~\eqref{eq:pref}.

These equality for the transformation of the series argument $z = (-1)^nn^nq^{n-1}/((n-1)^{n-1}p^n)$ to $z/(z-1)=\eta$ and
the equality for the prefactor term $-(q/p)(1-z)^{-a}$, are suggestive that we are on the right track.
The results suggest that there is indeed an identity analogous to Kummer's work, for higher hypergeometric series ${}_nF_{n-1}$ with $n\ge3$.
In the specific case of a series for the root of the trinomial eq.~\eqref{eq:tri}, we derived the coefficients of the relevant series.
(Actually, we derived a recurrence to calculate the series coefficients.)
It remains to generalize these results to higher hypergeometric series ${}_nF_{n-1}$ with arbitrary parameters, independent of any specific application.


\appendix  
\setcounter{equation}{0}
\section{Identities for hypergeometric functions}\label{sec:hyp}
Let $F(a,b;c;z)$ denote a ${}_2F_1$ hypergeometric series.
Kummer published the following identities:
\begin{subequations}
\begin{align}
\label{eq:Kum_a}
F(a,b;c;z) &= (1-z)^{-a}F\Bigl(a,c-b;c;\frac{z}{z-1}\Bigr) \,,
\\
\label{eq:Kum_b}
F(a,b;c;z) &= (1-z)^{-b}F\Bigl(c-a,b;c;\frac{z}{z-1}\Bigr) \,,
\\
\label{eq:Kum_c}
F(a,b;c;z) &= (1-z)^{c-a-b}F(c-a,c-b;c;z) \,.
\\
\label{eq:Kum_d}
F(a,b;c;z) &= \frac{\Gamma(c)\Gamma(c-a-b)}{\Gamma(c-a)\Gamma(c-b)}F(a,b;a+b-c+1;1-z)
\\ \nonumber
&\quad +(1-z)^{c-a-b}\frac{\Gamma(c)\Gamma(a+b-c)}{\Gamma(a)\Gamma(b)}F(c-a,c-b;c-a-b+1;1-z) \,.
\end{align}
\end{subequations}
We also note the following identities for the sum and difference of two binomial series:
\begin{subequations}
\begin{align}
\label{eq:hyp_sum}
(1+z)^{-2a} + (1-z)^{-2a} &= 2\,F\Bigl(a,a+\frac12;\frac12;z^2\Bigr) \,,
\\
\label{eq:hyp_diff}
(1+z)^{-2a} - (1-z)^{-2a} &= -4az\,F\Bigl(a+\frac12,a+1;\frac32;z^2\Bigr) \,.
\end{align}
\end{subequations}

\setcounter{equation}{0}
\section{Birkeland's series for the cubic}\label{sec:Birk_cubic}
Birkeland \cite{Birkeland_algebraic_1927} derived series solutions for the trinomial equation $x^n = gx^s +\beta$,
including solutions for the special case of the cubic.
We set $n=3$, $s=1$, $g=-p$ and $\beta=-q$ when quoting Birkeland's solutions for the cubic below.
The cubic equation is then $x^3+px+q=0$.
Birkeland derived series solutions for the two domains $|27q^2| < |4p^3|$ and $|27q^2| > |4p^3|$.

We begin with the domain $|27q^2| < |4p^3|$.
In terms of our notation, Birkeland's series for the roots are
\begin{subequations}
\label{eq:Bx123}
\begin{align}
\label{eq:Bx1}
x_1 &= \frac{q}{2p}\,F\Bigl(\frac13,\frac23;\frac32;-\frac{27q^2}{4p^3}\Bigr) 
-\sqrt{-p}\,F\Bigl(-\frac16,\frac16;\frac12;-\frac{27q^2}{4p^3}\Bigr) \,,
\\
\label{eq:Bx2}
x_2 &= \frac{q}{2p}\,F\Bigl(\frac13,\frac23;\frac32;-\frac{27q^2}{4p^3}\Bigr) 
+\sqrt{-p}\,F\Bigl(-\frac16,\frac16;\frac12;-\frac{27q^2}{4p^3}\Bigr) \,,
\\
\label{eq:Birk_x3}
x_3 &= -\frac{q}{p}\,F\Bigl(\frac13,\frac23;\frac32;-\frac{27q^2}{4p^3}\Bigr) \,.
\end{align}
\end{subequations}
We begin with $x_3$.
Applying Kummer's transformation in eq.~\eqref{eq:Kum_d} yields
\begin{equation}  
\begin{split}
  F\Bigl(\frac13,\frac23;\frac32;-\frac{27q^2}{4p^3}\Bigr)
&= \frac{\Gamma(\frac32)\Gamma(\frac12)}{\Gamma(\frac76)\Gamma(\frac56)}F\Bigl(\frac13,\frac23;\frac32;1+\frac{27q^2}{4p^3}\Bigr)
\\
&\quad +\Bigl(1+\frac{27q^2}{4p^3}\Bigr)^{1/2}\frac{\Gamma(\frac32)\Gamma(-\frac12)}{\Gamma(\frac13)\Gamma(\frac23)}F\Bigl(\frac76,\frac56;\frac32;1+\frac{27q^2}{4p^3}\Bigr)
\\
&= \frac32\,F\Bigl(\frac13,\frac23;\frac12;-\frac{D}{4p^3}\Bigr)
-\frac{\sqrt3}{2}\,\Bigl(-\frac{D}{4p^3}\Bigr)^{1/2}F\Bigl(\frac56,\frac76;\frac32;-\frac{D}{4p^3}\Bigr) \,.
\end{split}
\end{equation}
Then
\begin{equation}
\begin{split}
x_3 &= -\frac{3q}{2p}\,F\Bigl(\frac13,\frac23;\frac12;-\frac{D}{4p^3}\Bigr) 
+\frac{\sqrt{3}q}{2p}\,\Bigl(-\frac{D}{4p^3}\Bigr)^{1/2}\,F\Bigl(\frac56,\frac76;\frac32;-\frac{D}{4p^3}\Bigr) \,.
\end{split}
\end{equation}
Apply Kummer's transformation in eq.~\eqref{eq:Kum_c} to the second term to obtain
\begin{equation}
\label{eq:Kummer_x3p}
\begin{split}
x_3 &= -\frac{3q}{2p}\,F\Bigl(\frac13,\frac23;\frac12;-\frac{D}{4p^3}\Bigr) 
+\frac{\sqrt{D}}{6p}\,F\Bigl(\frac13,\frac23;\frac32;-\frac{D}{4p^3}\Bigr) \,.
\end{split}
\end{equation}
The series in eq.~\eqref{eq:Kummer_x3p} converge absolutely if and only if $|D| < 4|p|^3$.
This constraint requires $p<0$. This fact will be significant below.  

Next, suppose $q>0$ and let us process the root $x_1$.
Let us apply Kummer's transformation in eq.~\eqref{eq:Kum_d}.
Note that $1-(-27q^2)/(4p^3) = -D/(4p^3)$.
Then
\begin{subequations}
\begin{align}
F\Bigl(\frac13,\frac23;\frac32;-\frac{27q^2}{4p^3}\Bigr) &= \frac32\,F\Bigl(\frac13,\frac23;\frac12;-\frac{D}{4p^3}\Bigr)
-\frac{\sqrt3}{2}\sqrt{-\frac{D}{4p^3}}\,F\Bigl(\frac56,\frac76;\frac32;-\frac{D}{4p^3}\Bigr) \,,
\\
F\Bigl(-\frac16,\frac16;\frac12;-\frac{27q^2}{4p^3}\Bigr) &= \frac{\sqrt3}{2}\,F\Bigl(-\frac16,\frac16;\frac12;-\frac{D}{4p^3}\Bigr) 
+\frac16\sqrt{-\frac{D}{4p^3}}\,F\Bigl(\frac13,\frac23;\frac32;-\frac{D}{4p^3}\Bigr) \,.
\end{align}
\end{subequations}
The transformed expression for the root $x_1$ is
\begin{equation}  
\label{eq:process_x1}
\begin{split}
x_1 &= \frac{q}{2p}\,\biggl[\, \frac32\,F\Bigl(\frac13,\frac23;\frac12;-\frac{D}{4p^3}\Bigr)
-\frac{\sqrt3}{2}\sqrt{-\frac{D}{4p^3}}\,F\Bigl(\frac56,\frac76;\frac32;-\frac{D}{4p^3}\Bigr) \biggr]
\\
&\quad
-\sqrt{-p}\,\biggl[\, \frac{\sqrt3}{2}\,F\Bigl(-\frac16,\frac16;\frac12;-\frac{D}{4p^3}\Bigr) 
+\frac16\sqrt{-\frac{D}{4p^3}}\,F\Bigl(\frac13,\frac23;\frac32;-\frac{D}{4p^3}\Bigr) \biggr]
\\
&= \frac{3q}{4p}\,F\Bigl(\frac13,\frac23;\frac12;-\frac{D}{4p^3}\Bigr)
-\frac{\sqrt{-3p}}{2}\,F\Bigl(-\frac16,\frac16;\frac12;-\frac{D}{4p^3}\Bigr) 
\\
&\quad
-\sqrt{-\frac{D}{4p^3}}\biggl[ \frac{\sqrt3}{4}\,\frac{q}{p}\,F\Bigl(\frac56,\frac76;\frac32;-\frac{D}{4p^3}\Bigr) 
+\frac{\sqrt{-p}}{6}\,F\Bigl(\frac13,\frac23;\frac32;-\frac{D}{4p^3}\Bigr) \biggr] \,.
\end{split}
\end{equation}
The last two terms cancel.
Let us denote the sum of those two terms by $S$.
To prove $S=0$, we apply Kummer's transformation in eq.~\eqref{eq:Kum_c} to the second term:
\begin{equation}  
\begin{split}
S &= \frac{\sqrt3}{4}\,\frac{q}{p}\,F\Bigl(\frac56,\frac76;\frac32;-\frac{D}{4p^3}\Bigr) 
+\frac{\sqrt{-p}}{6}\,F\Bigl(\frac13,\frac23;\frac32;-\frac{D}{4p^3}\Bigr)
\\
&= \frac{\sqrt3}{4}\,\frac{q}{p}\,F\Bigl(\frac56,\frac76;\frac32;-\frac{D}{4p^3}\Bigr) 
  +\frac{\sqrt{-p}}{6}\Bigl(1+\frac{D}{4p^3}\Bigr)^{1/2} F\Bigl(\frac76,\frac56;\frac32;-\frac{D}{4p^3}\Bigr)
\\
&= \frac{\sqrt3}{4}\,\frac{q}{p}\,F\Bigl(\frac56,\frac76;\frac32;-\frac{D}{4p^3}\Bigr) 
  +\frac{\sqrt{-p}}{6}\Bigl(-\frac{27q^2}{4p^3}\Bigr)^{1/2} F\Bigl(\frac76,\frac56;\frac32;-\frac{D}{4p^3}\Bigr)
\\
&= \frac{\sqrt3}{4}\biggl(\frac{q}{p} +\frac{|q|}{|p|}\biggr) F\Bigl(\frac56,\frac76;\frac32;-\frac{D}{4p^3}\Bigr)
\\
&= 0 \,.
\end{split}
\end{equation}
This proves the cancellation, because $p<0$.
Hence
\begin{equation}  
x_1 = \frac{3q}{4p}\,F\Bigl(\frac13,\frac23;\frac12;-\frac{D}{4p^3}\Bigr)
-\frac{\sqrt{-3p}}{2}\,F\Bigl(-\frac16,\frac16;\frac12;-\frac{D}{4p^3}\Bigr) \,.
\end{equation}
We next apply Kummer's transformation in eq.~\eqref{eq:Kum_c} to the second term
(in the second line, $\sqrt{p^2}=|p|=-p$)
\begin{equation}  
\begin{split}
x_1 &= \frac{3q}{4p}\,F\Bigl(\frac13,\frac23;\frac12;-\frac{D}{4p^3}\Bigr)
-\frac{\sqrt{-3p}}{2}\,\Bigl(1+\frac{D}{4p^3}\Bigr)^{1/2}\,F\Bigl(\frac23,\frac13;\frac12;-\frac{D}{4p^3}\Bigr)
\\
&= \biggl(\frac{3q}{4p} -\frac12\sqrt{\frac{81q^2}{4p^2}}\biggr)\,F\Bigl(\frac13,\frac23;\frac12;-\frac{D}{4p^3}\Bigr)
\\
&= \biggl(\frac{3q}{4p} +\frac{9q}{4p}\biggr)\,F\Bigl(\frac13,\frac23;\frac12;-\frac{D}{4p^3}\Bigr)
\\
&= \frac{3q}{p}\,F\Bigl(\frac13,\frac23;\frac12;-\frac{D}{4p^3}\Bigr) \,.
\end{split}
\end{equation}
Correspondingly, if $q<0$, we process the root $x_2$ instead, to obtain the cancellation of terms obtained in eq.~\eqref{eq:process_x1}, to obtain a real root.

The three roots sum to zero, which yields the expression for the third root.
We denote the roots by $(r_1,r_2,r_3)$, where $r_1$ is always real.
The expressions for the three roots are
\begin{subequations}
\label{eq:B_disc_r123}
\begin{align}
r_1 &= \frac{3q}{p}\,F\Bigl(\frac13,\frac23;\frac12;-\frac{D}{4p^3}\Bigr) \,,
\\
r_2 &= -\frac{r_1}{2}
+\frac{\sqrt{D}}{6p}\,F\Bigl(\frac13,\frac23;\frac32;-\frac{D}{4p^3}\Bigr) \,,
\\
r_3 &= -\frac{r_1}{2}
-\frac{\sqrt{D}}{6p}\,F\Bigl(\frac13,\frac23;\frac32;-\frac{D}{4p^3}\Bigr) \,.
\end{align}
\end{subequations}
These expressions are the same as the roots in eqs.~\eqref{eq:r0prime} and \eqref{eq:r12prime}.

Next let us treat the domain $|27q^2| > |4p^3|$.
Let $\omega_3$ denote a primitive third root of unity.
Birkeland \cite{Birkeland_algebraic_1927} derived three series as follows (here $\ell=1,2,3$).
The expression below corrects some misprints in \cite{Birkeland_algebraic_1927}, and is
\begin{equation}
\label{eq:Broots_case2}  
\begin{split}
x_\ell &= -\frac{q}{|q|}\ \biggl[ \omega_3^\ell |q|^{1/3}\,{}_2F_1\Bigl(-\frac16,\frac13;\frac23;-\frac{4p^3}{27q^2}\Bigr)
    -\omega_3^{2\ell} \frac{p}{3|q|^{1/3}}\,{}_2F_1\Bigl(\frac16,\frac23;\frac43;-\frac{4p^3}{27q^2}\Bigr) \biggr] \,.
\end{split}
\end{equation}
Let us apply Kummer's transformation in eq.~\eqref{eq:Kum_a}.
First
\begin{equation}
\begin{split}
F\Bigl(-\frac16,\frac13;\frac23;-\frac{4p^3}{27q^2}\Bigr) 
&= \Bigl(1+\frac{4p^3}{27q^2}\Bigr)^{1/6}F\Bigl(-\frac16,\frac13;\frac23;-\frac{4p^3}{D}\Bigr) 
\\
&= \Bigl(-\frac{D}{27q^2}\Bigr)^{1/6}F\Bigl(-\frac16,\frac13;\frac23;-\frac{4p^3}{D}\Bigr) 
\\
&= \frac{(-D)^{1/6}}{\sqrt3\,|q|^{1/3}}\,F\Bigl(-\frac16,\frac13;\frac23;-\frac{4p^3}{D}\Bigr) \,.
\end{split}
\end{equation}
The above series converges absolutely if and only if $|D| > |4p^3|$.
This implies $D<0$.
Hence $(-D)^{1/6}$ is real in the above expressions.
Next
\begin{equation}
\begin{split}
F\Bigl(\frac16,\frac23;\frac43;-\frac{4p^3}{27q^2}\Bigr) 
&= \Bigl(1+\frac{4p^3}{27q^2}\Bigr)^{-1/6}F\Bigl(\frac16,\frac23;\frac43;-\frac{4p^3}{D}\Bigr) 
\\
&= \Bigl(-\frac{D}{27q^2}\Bigr)^{-1/6}F\Bigl(\frac16,\frac23;\frac43;-\frac{4p^3}{D}\Bigr) 
\\
&= \frac{\sqrt3|q|^{1/3}}{(-D)^{1/6}}F\Bigl(\frac16,\frac23;\frac43;-\frac{4p^3}{D}\Bigr) \,.
\end{split}
\end{equation}
Then
\begin{equation}
\begin{split}
x_\ell &= -\frac{1}{\sqrt3}\frac{q}{|q|}\ \biggl[ \omega_3^\ell (-D)^{1/6}\,F\Bigl(-\frac16,\frac13;\frac23;-\frac{4p^3}{D}\Bigr)
-\omega_3^{2\ell} \frac{p}{(-D)^{1/6}}\,F\Bigl(\frac16,\frac23;\frac43;-\frac{4p^3}{D}\Bigr) \biggr] \,.
\end{split}
\end{equation}
Then $x_1$ and $x_2$ are complex conjugates and $x_3$ is real.
Let us define sum and differences as follows
\begin{subequations}
\begin{align}
x_s &= -\frac{q}{\sqrt3|q|}\ \biggl[ (-D)^{1/6}\,F\Bigl(-\frac16,\frac13;\frac23;-\frac{4p^3}{D}\Bigr)
  + \frac{p}{(-D)^{1/6}}\,F\Bigl(\frac16,\frac23;\frac43;-\frac{4p^3}{D}\Bigr) \biggr] \,,
\\
x_d &= -\frac{q}{\sqrt{3}|q|}\ \biggl[ (-D)^{1/6}\,F\Bigl(-\frac16,\frac13;\frac23;-\frac{4p^3}{D}\Bigr)
  - \frac{p}{(-D)^{1/6}}\,F\Bigl(\frac16,\frac23;\frac43;-\frac{4p^3}{D}\Bigr) \biggr] \,.
\end{align}
\end{subequations}
Let us denote the roots by $(\rho_1,\rho_2,\rho_3)$, where $\rho_1$ is real.
The three roots are, in terms of hypergeometric series in the discriminant,
\begin{subequations}
\label{eq:B_disc_rho123}
\begin{align}
\rho_1 &= x_d \,,
\\
\rho_2 &= -\frac12\,x_d +i\frac{\sqrt{3}}{2}\,x_s \,,
\\
\rho_3 &= -\frac12\,x_d -i\frac{\sqrt{3}}{2}\,x_s \,.
\end{align}
\end{subequations}
Note the following.
\begin{enumerate}
\item
The root $\rho_1$ is real and $\rho_2$ and $\rho_3$ form a complex conjugate pair.
\item
In this domain, it is manifest from Birkeland's hypergeometric series in eq.~\eqref{eq:Broots_case2} that $x_3$ is real and $x_1$ and $x_2$ form a complex conjugate pair.
\item
The series in eq.~\eqref{eq:B_disc_rho123} converge absolutely if and only if $|D| > |4p^3|$.
This is not the same as the domain $|27q^2| > |4p^3|$ for Birkeland's series in eq.~\eqref{eq:Broots_case2}.
\end{enumerate}

\setcounter{equation}{0}
\section{Alternative discriminant series}\label{sec:dseries_alt}
We derive an alternative series for a root of the trinomial in eq.~\eqref{eq:tri} in powers of the discriminant, for $n>3$.
We assume $n$ is odd, so the discriminant has all negative coefficients, i.e., $D = -(n-1)^{n-1}p^n -n^nq^{n-1}$.
For brevity define $\xi = -D/((n-1)^{n-1}p^n)$.
Motivated by the series solution for the real root for a cubic,
we try a series of the following form, where $G_0$ is a power series 
\begin{equation}
r = \frac{sq}{p}\,G_0(\xi) \,.
\end{equation}
Here $s$ is a scale factor, to be determined below.
Let
\begin{equation}
G_0 = 1 + \sum_{k=1}^\infty \gamma_k \xi^k \,.
\end{equation}
Here $\gamma_0=1$ and $\{\gamma_k, k=1,2,\dots\}$ are a set of coefficients to be determined.
Substitute into eq.~\eqref{eq:tri} to obtain
\begin{equation}
\frac{s^nq^n}{p^n}G_0^n +sqG_0 + q = 0 \,.
\end{equation}
Retain terms in $G_0$ to the zeroth order only.
Then we obtain zeroth order term, say $T_0$:
\begin{equation}  
T_0 = \frac{s^nq^n}{p^n} +(s+1)q = \frac{q}{p^n}(s^nq^{n-1} +(s+1)p^n) \,.
\end{equation}
We demand that $T_0 \propto q\xi$, so that it can be absorbed/cancelled by the next-order term in $G_0$ of $O(\xi)$.
Hence $s$ must satisfy the constraint equation
\begin{equation}  
\frac{s^n}{s+1} = \frac{n^n}{(n-1)^{n-1}} \,.
\end{equation}
We know from the cubic trinomial that if $n=3$ then $s=3$.
If $n=5$, then $s\simeq 2.0632865$.
Observe that the determination of the value of $s$ itself requires the solving of a trinomial equation
\begin{equation}  
(n-1)^{n-1} s^n -n^n(s+1) = 0 \,.
\end{equation}
This is a weak point in the formalism, in that we must solve a specific version of eq.~\eqref{eq:tri} to proceed.
Anyway, assume we know the value of $s$.
Then 
\begin{equation}  
\begin{split}
T_0 &= \frac{s^nq^n}{p^n} +(s+1)q 
\\
&= q(s+1)\biggl[\frac{n^n}{(n-1)^{n-1}}\frac{q^{n-1}}{p^n} +1\biggr]
\\
&= \frac{q(s+1)}{(n-1)^{n-1}p^n}(n^n q^{n-1} +(n-1)^{n-1}p^n)
\\
&= q(s+1)\xi \,.
\end{split}
\end{equation}
Note also that
\begin{equation}
\begin{split}
\frac{q^n}{p^n} = \frac{n^nq^n}{n^np^n}
= q\frac{-D-(n-1)^{n-1}p^n}{n^np^n}
= q(\xi-1)\frac{(n-1)^{n-1}}{n^n} \,.
\end{split}
\end{equation}
Hence
\begin{equation}
\frac{s^nq^n}{p^n} = q(s+1)(\xi-1)\,.
\end{equation}
Also define a multinomial sum as follows.
Define $M_0=M_1=0$ and for $k\ge2$:
\begin{equation}
\label{eq:Mk}  
M_k = \sum_{j_1+2j_2+\dots+(k-1)j_{k-1}=k}\frac{n!}{j_1!j_2!\dots j_{k-1}!(n-j_1-j_2-\dots-j_{k-1})!} \gamma_1^{,j_1}\gamma_2^{j_2}\dots\gamma_{k-1}^{j_{k-1}} \,.
\end{equation}
Then expand terms to obtain
\begin{equation}
\begin{split}
0 &= \frac{s^nq^n}{p^n}\biggl(1 +\sum_{k=1}^\infty \gamma_k\xi^k\biggr)^n +sq\biggl(1+\sum_{k=1}^\infty \gamma_k\xi^k\biggr) +q
\\
&= \frac{s^nq^n}{p^n}\biggl[1 +n\gamma_1\xi +(n\gamma_2+M_2)\xi^2 +(n\gamma_3+M_3)\xi^3 +\dots\biggr]
\\
&\qquad +sq(\gamma_1\xi+\gamma_2\xi^2+\dots) + q(s+1)
\\
&= q(s+1)(\xi-1)\biggl[n\gamma_1\xi +(n\gamma_2+M_2)\xi^2 +(n\gamma_3+M_3)\xi^3 +\dots\biggr]
\\
&\qquad
+sq(\gamma_1\xi+\gamma_2\xi^2+\dots) +q(s+1)\xi \,.
\end{split}
\end{equation}
We equate powers of $\xi$. First
\begin{equation}
\begin{split}
0 &= -(s+1)n\gamma_1 +s\gamma_1 +(s+1)
\\
\Rightarrow\qquad (ns+n-s)\gamma_1 &= s+1
\\
\Rightarrow\qquad \gamma_1 &= \frac{s+1}{(n-1)(s+1)+1} \,.
\end{split}
\end{equation}
For higher powers of $\xi$, we obtain a recurrence for $\gamma_k$:
\begin{equation}
\begin{split}
0 &= (s+1)(n\gamma_{k-1}+M_{k-1}) -(s+1)(n\gamma_k+M_k) +s\gamma_k
\\
\Rightarrow\qquad (ns+n-s)\gamma_k &= (s+1)(n\gamma_{k-1} +M_{k-1} -M_k) 
\\
\Rightarrow\qquad \gamma_k &= \frac{s+1}{(n-1)(s+1)+1}(n\gamma_{k-1} +M_{k-1} -M_k) \,.
\end{split}
\end{equation}
The above scheme has been tested numerically for $n=3$, $5$ and $7$ and it works.

We postulate that the series $G_0$ is a ${}_{n-1}F_{n-2}$ hypergeometric series and the root $r$ has the form
\begin{equation}  
r = \frac{sq}{p}{}_{n-1}F_{n-2}(a_n,2a_n,\dots,(n-1)a_n;c_n,2c_n,\dots(n-2)c_n; \xi) \,.
\end{equation}
Here $a_n$ and $c_n$ are parameters to be determined.
For $n=3$, we know that $s=3$.
Working out the coefficients $\gamma_k$ yields that the root is given by the hypergeometric series
\begin{equation}  
  r = \frac{3q}{p}F\Bigl(\frac13,\frac23;\frac12;-\frac{D}{4p^3}\Bigr)\,.
\end{equation}
This equals the root in eq.~\eqref{eq:r0prime}.
For $n>3$, note that
\begin{subequations}  
\begin{align}  
\gamma_1 &= \frac{s+1}{(n-1)(s+1)+1} \,,
\\
\gamma_2 &= \frac{\gamma_1}{2}\biggl[1+\frac{s}{(n-1)(s+1)+1}\biggr] \biggl[1 + \frac{1}{(n-1)(s+1)+1} \biggr] \,.
\end{align}  
\end{subequations}
We fit values for $a_n$ and $c_n$ using the above expressions for $\gamma_1$ and $\gamma_2$. 
\begin{enumerate}
\item
For $n=5$, we find $s_5\simeq 2.06328648929924$ and
\begin{equation}  
a_5 \simeq 0.06488215078404166 \,, \qquad
c_5 \simeq  0.06743687485964298 \,,
\end{equation}  
and
\begin{equation}  
r = \frac{s_5q}{p}{}_4F_3\Bigl(a_5,2a_5,3a_5,4a_5;c_5,2c_5,3c_5;-\frac{D}{4^4p^5}\Bigr) \,.
\end{equation}
\item
For $n=7$, we find $s_5\simeq 1.7404309861632174957$ and
\begin{equation}  
a_7 \simeq 0.02844486139350898 \,, \qquad
c_7 \simeq 0.0289201297905275 \,,
\end{equation}  
and
\begin{equation}  
r = \frac{s_7q}{p}{}_6F_5(a_7,2a_7,3a_7,4a_7,5a_7,6a_7;c_7,2c_7,3c_7,4c_7,5c_7;-\frac{D}{6^6p^5}\Bigr) \,.
\end{equation}
\end{enumerate}

The above scheme does not work for even $n$.
For even $n$, the constraint equation on $s$ is
\begin{equation}  
\frac{s^n}{s+1} = -\frac{n^n}{(n-1)^{n-1}} \,.
\end{equation}
The only real solution is a negative number $s = -n/(n-1)$.
However, the denominator $(n-1)(s+1)+1$ equals zero if $s=-n/(n-1)$, leading to a division by zero in the formalism.
It is left for future research to find a partner series in the discriminant for even $n$.



\begin{thebibliography}{0}

\bibitem{Birkeland_algebraic_1927}
R.~Birkeland, {\"U}ber die Aufl{\"o}sung algebraischer Gleichungen durch hypergeometrische Funktionen.
\textit{Mathematische Zeitschrift}, {\bf 26} (1927) 566--578.

\bibitem{Euler_1779}
L.~Euler, ``De serie Lambertina plurimisque eius insignibus proprietatibus,''
{\it Acta Academiae Scientarum Imperialis Petropolitinae}, {\bf 2} (1783) (original date 1779) 29--51.
Reprinted in L.~Euler, {\it Opera Omnia, Series Prima} in
``Commentationes Algebraicae,'' Teubner, Leipzig, Germany, {\bf 6} (1921) 350--369.

\bibitem{Lambert_1758}
J.~Lambert, ``Observationes variae in mathesin puram,'' {\it Acta Helvetica}, {\bf 3} (1758) 128--168.

\bibitem{Lambert_1770}
J.~Lambert, ``Observations analytiques,'' {\it Nouveaux M{\'e}moires de l'Acad{\'e}mie Royale des Sciences et Belles-Lettres},
Berlin, (1770) pp.~225--244.

\bibitem{Mane_16_1}
S.~R.~Mane, Multiparameter Fuss-Catalan numbers with application to algebraic equations.
\textit{arXiv:2410.07922 [math.CO]}, (2016).

\bibitem{Zucker}
I.~J.~Zucker, The cubic equation -- a new look at the irreducible case.
\textit{The Mathematical Gazette}, {\bf 92} (2008) 264--268.
\end{thebibliography}
\end{document}